\documentclass[aos,MSNbibl,citesort,dvips]{arximspdf}
\usepackage{mathbh}
\usepackage{graphicx}

\doi{10.1214/11-AOS907}
\volume{39}
\issue{5}
\pubyear{2011}
\firstpage{2477}
\lastpage{2501}

\makeatletter

\newtheorem{lemma}{Lemma}[section]
\newtheorem{theorem}{Theorem}[section]
\newtheorem{proposition}{Proposition}[section]
\newproclaim{definition}{Definition}[section]
\newproclaim{problem}{Problem}

\newcommand{\Risk}{\operatorname{Risk}}

\newcommand{\rd}{\mathrm{d}}
\renewcommand{\hat}{\widehat}
\renewcommand{\tilde}{\widetilde}
\newcommand{\vf}{\varphi}
\newcommand{\e}{\epsilon}

\newcommand{\cF}{\mathcal{F}}
\newcommand{\cX}{\mathcal{X}}

\newcommand{\bE}{\mathrm{E}}
\newcommand{\bL}{{\mathbb L}}
\newcommand{\bP}{\mathrm{P}}
\newcommand{\bR}{{\mathbb R}}
\newcommand{\bS}{{\mathbb S}}

\newcommand{\Si}{\operatorname{Si}}
\renewcommand{\P}{\mathcal{P}}
\newcommand{\on}{\frac{1}{n}}
\newcommand{\ind}{\mathbh{1}} 
\newcommand{\eep}{\varepsilon}
\newcommand{\wt}[1]{\widetilde{#1}}
\newcommand{\half}{\frac{1}{2}}
\newcommand{\rE}{\mathrm{E}}
\newcommand{\rR}{\mathbb R}
\newcommand{\rP}{\mathrm{P}}
\renewcommand{\kappa}{\varkappa}

\makeatother

\begin{document}
\begin{frontmatter}

\title{On deconvolution of distribution functions}
\runtitle{On deconvolution of distribution functions}

\begin{aug}
\author[A]{\fnms{I.} \snm{Dattner}\thanksref{t1}\ead[label=e1]{idattner@stat.haifa.ac.il}},
\author[A]{\fnms{A.} \snm{Goldenshluger}\corref{}\thanksref{t1}\ead[label=e2]{goldensh@stat.haifa.ac.il}}
\and
\author[B]{\fnms{A.} \snm{Juditsky}\ead[label=e3]{anatoli.iouditski@imag.fr}}
\runauthor{I. Dattner, A. Goldenshluger and A. Juditsky}
\affiliation{University of Haifa, University of Haifa
and Universit\'e Grenoble I}
\address[A]{I. Dattner\\
A. Goldenshluger\\
Department of Statistics\\
University of Haifa\\
31905 Haifa\\
Israel\\
\printead{e1}\\
\hphantom{E-mail: }\printead*{e2}}
\address[B]{A. Juditsky\\
LMC, B. P. 53\\
Universit\'e Grenoble I\\
38041 Grenoble Cedex 9\\
France\\
\printead{e3}} 
\end{aug}

\thankstext{t1}{Supported by BSF Grant 2006075.}

\received{\smonth{6} \syear{2010}}
\revised{\smonth{5} \syear{2011}}

%
\begin{abstract}
The subject of this paper is the problem of nonparametric estimation of
a continuous distribution function from observations with measurement
errors. We study minimax complexity of this problem when unknown
distribution has a density belonging to the Sobolev class, and the
error density is ordinary smooth. We develop rate optimal estimators
based on direct inversion of empirical characteristic function. We also
derive minimax affine estimators of the distribution function which are
given by an explicit convex optimization
problem. 
Adaptive versions of these estimators are proposed, and 
some numerical results demonstrating good practical behavior of the
developed procedures are presented.
\end{abstract}

%

\end{frontmatter}

\section{Introduction}
In this paper we study the problem of estimating
a distribution function in the presence of measurement errors.

Let $X_1,\ldots, X_n$ be a sequence of independent, identically
distributed random variables with common distribution $F$. Suppose that we
observe random variables $Y_1, \ldots, Y_n$ given by
%
%
\begin{equation}\label{eqonesample-model}
Y_j=X_j+\zeta_j,\qquad j=1,\ldots,n,
\end{equation}
where $\zeta_j$ are i.i.d. random variables, independent of $X_j$'s with
the density $f_\zeta$ w.r.t. the Lebesgue measure on the real line.
The objective is to estimate the value $F(t_0)$ of the
distribution function $F$ of $X$ at a given point $t_0\in\bR$
from the observations $Y^n=(Y_1, \ldots, Y_n)$.

By an estimator we mean any measurable function
$\tilde{F}=\tilde{F}(Y^n)$ of the observations $Y^n$. We adopt the
minimax approach for measuring estimation accuracy.\vspace*{1pt} Let
$\cF$ be a given family of probability distributions on $\rR$. Given an
estimator $\tilde{F}$ of $F(t_0)$, we consider two types of maximal
over $\cF$ risks:
\begin{itemize}
\item quadratic risk,
\[
\Risk_2[\tilde{F}; \cF]:= \sup_{F\in\cF}
\{\rE|\tilde{F}-F(t_0)|^2\}^{1/2}.
\]
\item$\epsilon$-risk: given a tolerance level $\epsilon\in(0,1/2)$
we define
\[
\Risk_{\epsilon}[\tilde{F};\cF] := \min\Bigl\{\delta\dvtx\sup
_{F\in\cF} \bP
[|\tilde{F}-F(t_0)|>\delta]\leq\epsilon\Bigr\}.
\]
\end{itemize}
An estimator $\tilde{F}{}^*$ is said to be \textit{rate optimal} or
\textit{optimal in order}
with respect to $\Risk$
if
\[
\Risk[\tilde{F}^*; \cF] \leq C \inf_{\tilde{F}} \Risk[\tilde{F};
\cF],
\]
where $\inf$ is taken over all possible estimators of $F(t_0)$,
and $C<\infty$ is independent of $n$.
We will be particularly interested in the classes of distributions
having density with respect to the
Lebesgue measure on the real line.

The outlined problem is closely related to the density deconvolution
problem that has been extensively studied in the literature; see, for
example,
\cite{zhang,stefanski,Fan,butucea-tsy12,butucea-tsy,johannes,meister}
and references therein. In these works the minimax rates of convergence
have been derived under different assumptions on the error density and
on the smoothness of the density to be estimated. Depending on the tail
behavior of the characteristic function $\hat{f}_\zeta$ of $\zeta$ the
following two cases are usually distinguished:
\begin{longlist}[(ii)]
\item[(i)] \textit{ordinary smooth} errors, when the tails
of $\hat{f}_\zeta$ are polynomial, that is,
\[
|\hat{f}_\zeta(\omega)|\asymp|\omega|^{-\beta},\qquad|\omega|\to
\infty,
\]
for
some $\beta>0$;
\item[(ii)] \textit{supersmooth} errors, when the tails are exponential,
that is,
\[
|\hat{f}_\zeta(\omega)|\asymp\exp\{-c|\omega|^\beta\},\qquad
|\omega|\to\infty,
\]
for some $c>0$ and $\beta>0$.
\end{longlist}
The afore cited papers derive minimax rates of convergence
for different functional classes under ordinary smooth and supersmooth
errors.

In contrast to existence of the voluminous literature
on density deconvolution, the problem of deconvolution of the distribution
function $F$
has attracted much less attention and has been studied in very few
papers (see
\cite{meister}, Section 2.7.2,
for a recent review of corresponding contributions).
A consistent estimator of a distribution function
from observations with additive Gaussian measurement errors
was developed by \cite{gaffey}.
A ``plug-in'' estimator based on integration of the density
estimator in the density deconvolution problem has been studied under moment
conditions on $F$ in \cite{zhang}.
The paper \cite{Fan} also considered the estimator based on
integration of the
density deconvolution estimator. It was shown there that
under a tail condition on $F$ the
estimator achieves optimal rates
of convergence provided that the errors are supersmooth.
For the case of ordinary smooth errors there is a gap between the upper
and lower bounds reported in \cite{Fan} which leaves open the question
of constructing optimal estimators.
More recently, some
minimax rates of estimation of distribution functions
in models with measurement errors were 
reported
in \cite{hall}. Note also that \cite{butucea-comte} considered a general
problem of optimal and adaptive
estimation of linear functionals
$\ell(f)=\int_{-\infty}^{\infty}
\phi(t) f(t) \,\rd t$ in the model (\ref{eqonesample-model}). However, their
results hold only for representative $\phi\in\bL_1(\bR)$ which is
clearly not
the
case in the problem of recovery of distribution function.

The objective of this paper is to develop optimal methods of
minimax deconvolution of distribution functions and
to answer several questions raised
by known results on this problem:
Is a smoothness assumption alone on $F$ sufficient
in order to secure
minimax rates of estimation of the sort $O(n^{-\gamma})$ for $\gamma
>0$ in the
case of
ordinary smooth errors?
Do we need tail or moment conditions on $F$?
%

Our contribution is two-fold.
First, we characterize the minimax rates of convergence
in the case when the unknown distribution belongs to a
Sobolev ball, and the observation errors are ordinary smooth.
The rates of convergence depend crucially on
the relation between the \textit{smoothness index} $\alpha$ of the
Sobolev ball
and the parameter $\beta$ [the rate at which
the characteristic function of errors tends to zero; see (i) above].
In contrast to the density deconvolution problem, it turns out that
there are
different regions in the $(\alpha, \beta)$-plane
where different rates of convergence are attained.
We show that in some regions of the $(\alpha,\beta)$-plane
the minimax rates of convergence are attained by a linear
estimator, which is based
on direct inversion of the
distribution function from the corresponding characteristic function;
cf. \cite{hall}.
It is worth noting that we
do not require any additional tail or moment conditions on the
unknown
distribution.
In the case when the parameters
of the regularity class of the
distribution $F$ are unknown,
we also construct an adaptive estimator
based on Lepski's adaptation scheme \cite{lepski}.
The $\epsilon$-risk of this estimator is within a $\ln\ln n$-factor
of the minimax $\epsilon$-risk.

Second, using recent results on
estimating
linear functionals developed in \cite{JudNem},
we propose minimax and
adaptive affine estimators of the cumulative distribution function
for a discrete distribution deconvolution problem; see
also \cite{DonLiuI,DonLiuII,DonLiuIII,CaiLow2003,Donoho} for the general
theory of
affine estimation. These estimators can be applied
to the original deconvolution problem provided that it
can be
efficiently discretized. By efficient discretization we mean that:
\begin{longlist}[(2)]
\item
[(1)] the support of the distributions of $X$ ($Y$) can be ``compactified''
[one can point out a compact subset of $\bR$ such that the probability
of $X$ ($Y$) being outside this set is ``small'']
and binned into small intervals;
\item[(2)] the class $\cX$ of discrete distributions,
obtained by the corresponding finite-dimensional cross-section
of the class $\cF$ of continuous distributions
is a computationally tractable convex closed set.%
\setcounter{footnote}{1}\footnote{%
Roughly speaking, a computationally tractable
set can be interpreted as a set given by a finite system of inequalities
$p_i(x)\leq0$, $i=1,\ldots,m$, where $p_i$ are convex polynomials; see,
for example,
\cite{BN}, Chapter 4.}
\end{longlist}
Under these conditions one can efficiently implement
the minimax affine estimator
for $F$ based on the approach proposed in \cite{JudNem}.
This estimator is rate minimax with respect to $\Risk_\e$ (within
a factor $\approx2$ for small $\epsilon$) \textit{whatever are the noise
distribution and a convex and closed class $\cX$}.

We describe construction of the minimax
affine estimator of $F$
when the class $\cX$ is known and provide an
adaptive version of the estimation procedure
when the available information allows us to construct an embedded
family of
classes.

The rest of the paper is structured as follows. We present our results on
estimation over the
Sobolev classes in Section \ref{secminimax}.
Section \ref{secminimax-linear}
deals with minimax and adaptive affine estimation.
Section \ref{secnumeric} presents a numerical study of proposed adaptive
estimators and discusses
their relative merits.
Proofs of all results are given in the supplementary
article \cite{DGJ-suppl}.
%

%
\section{Estimation over Sobolev classes}\label{secminimax}
\subsection{Notation}
We denote by $f_Y$ and $f_\zeta$ the densities
of random variables $Y$ and $\zeta$; with certain
abuse of notation we simply denote by $f$ the density of unknown
distribution of $X$.

Let $g$ be a function on $\bR$; we denote by $\hat{g}$ the Fourier
transform of $g$,
\[
\hat{g}(\omega)=\int_{-\infty}^\infty g(x) e^{i\omega x}\,\rd x,\qquad
\omega\in\bR.
\]
We consider the classes of absolutely continuous distributions.
\begin{definition}\label{21}
Let $\alpha> -1/2$, $L>0$. We say that $F$ belongs to the class $\cF
_\alpha(L)$
if it has a density $f$ with respect to the Lebesgue measure on $\rR$, and
\[
\frac{1}{2\pi}
\int_{-\infty}^\infty|\hat{f}(\omega)|^2 (1+\omega^2)^{\alpha}
\,\rd\omega\leq L^2.
\]
\end{definition}

The set $\cF_\alpha(L)$ with $\alpha>-1/2$
contains absolutely continuous distributions.
If
$\alpha>1/2$, then the distributions $F$ from $\cF_\alpha(L)$ have
bounded continuous densities. Usually
$\cF_\alpha(L)$ is referred to
as the Sobolev class.

We use extensively
the following inversion formula: for a continuous distribution
$F$ one has
%
%
\begin{equation}\label{eqinversion}
F(x)=\frac{1}{2}-\frac{1}{\pi}\int_0^{\infty} \omega^{-1}
\Im\{e^{-i\omega x}\hat{f}(\omega)\} \,\rd\omega,\qquad x\in\bR,
\end{equation}
where $\Im\{\cdot\}$ stands for the imaginary part, and the above
integral is
interpreted as an improper Riemann integral
$
\lim_{T\rightarrow\infty}\int_{1/T}^T
\omega^{-1}\Im\{e^{-i\omega x}\hat{f}(\omega)\}\,\rd\omega.
$
For the proof of (\ref{eqinversion}) see \cite{gurland,gil-pelaez} and
\cite{kendall}, Section 4.3.

Throughout this section
we assume that the error characteristic function
does not vanish:
\[
|\hat{f}_\zeta(\omega)|\ne0\qquad\forall\omega\in\bR.
\]
This is a standard assumption in deconvolution problems.

\subsection{Minimax rates of estimation}\label{secminmax1}
In model (\ref{eqonesample-model}) we have
$\hat{f}(\omega)=\hat{f}_Y(\omega)/\hat{f}_\zeta(\omega)$,
and $\hat{f}_Y(\omega)$ can be easily estimated by the empirical
characteristic function of the observations
$Y$.
This motivates the following construction:
for $\lambda>0$ we define the estimator $\tilde{F}_{\lambda}$ of $F(t_0)$
by
%
%
\begin{equation}\label{eqcdfEstimator1}
\tilde{F}_{\lambda}=\frac{1}{2}-
\frac{1}{n}\sum_{j=1}^n \frac{1}{\pi}\int_0^\lambda\frac
{1}{\omega}
\Im\biggl\{
\frac{e^{i\omega(Y_j-t_0)}}{\hat{f}_\zeta(\omega)}
\biggr\}\,\rd\omega.
\end{equation}
Here $\lambda$ is the design parameter to be specified.
Note that if the density $f_\zeta$ is symmetric around the origin, then
$\hat{f}_\zeta$ is real, and the estimator $\tilde{F}_\lambda(t_0)$
takes the form (cf. \cite{hall})
\[
\tilde{F}_{\lambda}=\frac{1}{2}-
\frac{1}{n}\sum_{j=1}^n \frac{1}{\pi}\int_0^{\lambda}\frac{\sin
\{\omega(Y_j-t_0)\}}
{\hat{f}_\zeta(\omega)\omega}\,\rd\omega.
\]
Note that $\tilde{F}_\lambda$ may be truncated
to the interval $[0,1]$; obviously, the risk of such a ``projected''
estimator is
smaller than that of $\tilde{F}_\lambda$.

Our current goal is to establish an upper bound on the risk of the
estimator $\tilde{F}_\lambda$ over the classes $\cF_\alpha(L)$.
We need the following assumptions on the distribution of the measurement
errors $\zeta_i$:

\renewcommand\thelonglist{(E1)}
\renewcommand\labellonglist{\thelonglist}
\begin{longlist}
\item\label{assE1}
There exist real numbers
$\beta> 0$, $c_{\zeta}>0$ and $C_\zeta>0$ such that
\[
c_\zeta(1+ \omega^2)^{-\beta/2} \leq
|\hat{f}_\zeta(\omega)| \leq C_{\zeta}(1+ \omega^2)^{-\beta
/2}\qquad \forall\omega\in\bR.
\]
\end{longlist}
\renewcommand\thelonglist{(E2)}
\renewcommand\labellonglist{\thelonglist}
\begin{longlist}
\item\label{assE2}
There exist positive real numbers
$\omega_0$, $b_\zeta$ and $\tau$ such that
\[
|\hat{f}_\zeta(\omega)| \geq1- b_\zeta|\omega|^\tau\qquad
\forall|\omega|\leq
\omega_0.
\]
%
\end{longlist}

%

Assumption \ref{assE1} characterizes the case of the
\textit{ordinary smooth errors}.
Assumption \ref{assE2}
describes the local behavior of $\hat{f}_\zeta$ near the origin. It is
well known that
for any distribution of a nondegenerate random variable there exist positive
constants $b$ and $\delta$ such that
$|\hat{f}(\omega)|\leq1- b|\omega|^2$ for all $|\omega|\leq\delta$
(see, e.g., \cite{petrov}, Lemma 1.5). Thus in \ref{assE2} we have
$\tau\in(0, 2]$.
Typical examples of distributions satisfying \ref{assE1} and \ref
{assE2} are the Laplace
and Gamma distributions. For example, for the Laplace distribution
\ref{assE1} holds with $\beta=2$, and \ref{assE2} holds with $\tau=2$.
The Gamma distribution
provides
an example of the distribution satisfying \ref{assE1}
with $\beta>0$ being the
shape parameter of the distribution.

\begin{figure}

\includegraphics{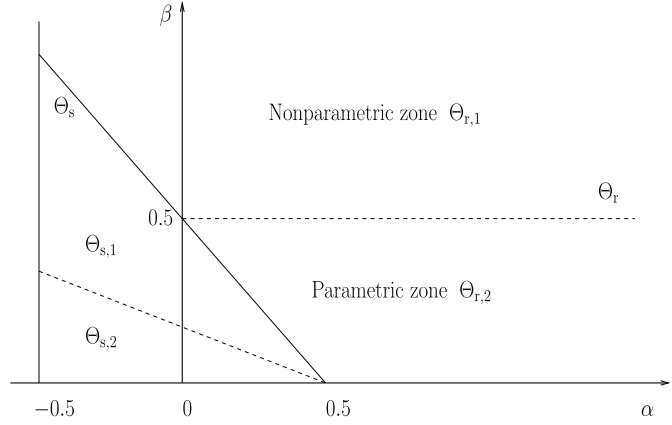}

\caption{Division of the parameter set for $(\alpha, \beta)$.}
\label{fig1}
\end{figure}

As we will see in the sequel, the rates of convergence of the risks
$\Risk_2[\tilde{F}_\lambda$; $\cF_\alpha(L)]$ and $\Risk_{\epsilon
}[\tilde{F}_\lambda; \cF_\alpha(L)]$
are mainly determined by the relationship between parameters
$\alpha$ and $\beta$.
Consider the following two subsets
of the parameter set
$\Theta:=\{(\alpha, \beta)\dvtx\alpha>- 1/2, \beta>0\}$
for the pair $(\alpha, \beta)$:
\[
\Theta_{\mathrm r}:=\{(\alpha, \beta)\in\Theta\dvtx\alpha+\beta
>1/2\},\qquad
\Theta_{\mathrm s}:=\{(\alpha,\beta)\in\Theta\dvtx\alpha+\beta
<1/2\}.
\]
If
$(\alpha,\beta)\in\Theta_{\mathrm s}$, then necessarily
$\hat{f}_\zeta\notin\bL_1(\bR)$; in addition, because
$\alpha<1/2$,
the density $f$ can be discontinuous.
That is why we will refer to $\Theta_{\mathrm s}$
as the \textit{singular zone}, while
the subset $\Theta_{\mathrm r}$ will be called the \textit{regular zone}.
We denote by $\Theta_{\mathrm b}$ the \textit{border zone} between
$\Theta_{\mathrm r}$ and $\Theta_{\mathrm s}$:
\[
\Theta_{\mathrm b}:=\{(\alpha,\beta)\in\Theta\dvtx\alpha+\beta=1/2\}.
\]
Division of the parameter set $\Theta$
into zones $\Theta_{\mathrm r}$, $\Theta_{\mathrm s}$ and $\Theta
_{\mathrm b}$
is displayed
in Figure~\ref{fig1}.
%
%
The figure also shows the sub-regions
$\Theta_{{\mathrm r},i}$ and $\Theta_{{\mathrm s},i}$, $i=1,2$,
that are defined by the following formulas:
\begin{eqnarray*}
\Theta_{{\mathrm r},1}&:=&\{(\alpha,\beta)\in\Theta_{\mathrm r}\dvtx\beta
>1/2\}
,\qquad
\Theta_{{\mathrm r}, 2}:=\{(\alpha,\beta)\in\Theta_{\mathrm r}\dvtx
\beta
<1/2\},
\\
\Theta_{{\mathrm s},1}&:=&\{(\alpha, \beta)\in\Theta_{\mathrm s}\dvtx
\alpha
+3\beta\geq1/2\},\qquad
\Theta_{{\mathrm s},2}:=\{(\alpha,\beta)\in\Theta_{\mathrm s}\dvtx
\alpha
+3\beta<1/2\}.
\end{eqnarray*}

The next
two theorems present bounds on the risks in the regular zone:
Theorem \ref{thupper-bound} states upper bounds
on the risks of $\tilde{F}_\lambda$, while Theorem \ref{thlower-1}
contains the corresponding lower bounds on the minimax risks.

For $z\geq1$ define
\[
\lambda(z)=z^{1/[2\alpha+ (2\beta\vee
1)]} ,\qquad
\psi(z)=\cases{
z^{-(2\alpha+1)/(4\alpha+4\beta)}, &\quad$\beta>1/2$, \vspace*{2pt}\cr
\sqrt{\ln z/z}, &\quad$\beta=1/2$,\vspace*{2pt}\cr
1/\sqrt{z}, &\quad$\beta\in(0, 1/2)$.}
\]

\begin{theorem}\label{thupper-bound}
Let\vspace*{1pt} assumptions \ref{assE1} and \ref{assE2} hold, and suppose that
$(\alpha, \beta)\in\Theta_{\mathrm r}$.
If
$\tilde{F}_{\lambda_\star}$
is estimator (\ref{eqcdfEstimator1}) associated with
$\lambda_\star=C_1(\alpha, L) \lambda(n)$, then
for all $t_0\in\bR$ and large enough $n$,
\[
\Risk_2[\tilde{F}_{\lambda_\star}; \cF_\alpha(L)] \leq
\psi_n(\alpha, L)
:= C_2(\alpha, L)\psi(n).
\]
In addition, if
$\lambda_\star=C_1(\alpha, L) \lambda(n/\ln[2\epsilon^{-1}])$, then
for all $t_0\in\bR$ and large enough~$n$,
\[
\Risk_{\epsilon}[\tilde{F}_{\lambda_\star}; \cF_\alpha(L)]
\leq
\psi_{n,\epsilon}(\alpha, L) := C_3(\alpha, L) \psi(n/\ln
[2\epsilon^{-1}]),
\]
provided that
$\epsilon\geq2\exp\{-C_4(\alpha, L) n\}$.
The constants
$C_i$, $i=1,\ldots, 4$, are specified in the proof of the theorem
(see \textup{(A.15)--(A.22)} in \cite{DGJ-suppl}).
\end{theorem}

Theorem \ref{thupper-bound} shows that if $(\alpha, \beta)$ is in
the regular zone
$\Theta_{\mathrm r}$
and $\beta\in(0,1/2)$, then the estimator $\tilde{F}_{\lambda_\star
}$ attains the parametric rate
of convergence. In the case $\beta=1/2$ this rate is within a
logarithmic factor
of the parametric rate. The natural question is
if the estimator $\tilde{F}_{\lambda_\star}$ is rate optimal
whenever $\beta>1/2$, and
$(\alpha, \beta)\in\Theta_{\mathrm r}$. The answer is provided by
Theorem \ref{thlower-1}.

We need the following assumption.
\renewcommand\thelonglist{(E3)}
\renewcommand\labellonglist{\thelonglist}
\begin{longlist}
\item\label{assE3}
The characteristic function $\hat{f}_\zeta$ is
twice differentiable, and there exist real numbers $\beta>1/2$,
$C_\zeta>0$ and $\omega_*>0$ such that
\[
(1+\omega^2)^{\beta/2} \max_{j=0,1,2} \bigl\{
\bigl|\hat{f}^{\,(j)}_\zeta(\omega)\bigr|\bigr\} \leq C_\zeta\qquad
\forall|\omega| \geq\omega_*.
\]
\end{longlist}
%
%
Assumption \ref{assE3} is rather standard in derivations of lower bounds
for deconvolution problems.
This assumption should be compared to
condition (G3) in \cite{Fan}; it is assumed there that for
$j=0,1,2$ one has
$|\hat{f}^{\,(j)}_\zeta(\omega)| |\omega|^{\beta+j}\leq C_\zeta$
as $|\omega|\to\infty$. Note that \ref{assE3} is a weaker assumption.
\begin{theorem}\label{thlower-1}
Let assumption \ref{assE3} hold.
Suppose that the class $\cF_\alpha(L)$ is such that
$L^2\geq
\pi^{-1}2^{1+(\alpha-1)_+}\Gamma(2\alpha+1)$ and $\alpha>1/2$. Then
there exist constants $c_1$ and $c_2$
depending
on $\alpha$, $\beta$ and $f_\zeta$ only such that, for all $n$ large enough,
\begin{eqnarray*}
\inf_{\tilde{F}}
\Risk_2[\tilde{F}; \cF_\alpha(L)] &\geq& c_1L^{(2\beta
-1)/(2\alpha+2\beta)}\phi_n,
\\
\inf_{\tilde{F}}\Risk_{\epsilon}[\tilde{F}; \cF_\alpha(L)]
&\geq& c_2 L^{(2\beta-1)/(2\alpha+2\beta)}
\phi_{n,\epsilon},
\end{eqnarray*}
where
$\phi_n:=\phi(n)$, $\phi_{n,\epsilon}:=\phi(n/\ln\epsilon^{-1})$,
$\phi(z):= z^{-(2\alpha+1)/(4\alpha+4\beta)}$,
and $\inf$ is taken over all possible estimators of $F(t_0)$.
\end{theorem}

The results of Theorems \ref{thupper-bound} and \ref{thlower-1} deal
with the regular zone. While we do not present the lower bound for the
case of
$\alpha\le1/2$ we expect that the bounds of Theorem \ref{thlower-1}
hold for
the whole regular zone.

It is important to realize that
the risks of $\tilde{F}_\lambda$
converge to zero for all $(\alpha, \beta)\in\Theta$, and, in
particular, for
$(\alpha, \beta)\in\Theta_{\mathrm s}$ and $(\alpha, \beta)\in\Theta
_{\mathrm b}$.
The next\vspace*{1pt} statement establishes upper bounds on $\Risk_2[\tilde
{F}_\lambda; \cF_\alpha(L)]$
in the singular and border zones, $\Theta_{\mathrm s}$ and $\Theta
_{\mathrm b}$.
\begin{theorem} \label{thupper-bound-2}
Let assumptions \ref{assE1} and \ref{assE2} hold.
If
$\tilde{F}_{\lambda_\star}$
is the estimator (\ref{eqcdfEstimator1}) associated with
$\lambda_\star=C_1(\alpha, L) \lambda(n)$, then
for all $t_0\in\bR$ and large enough $n$
\[
\Risk_2[\tilde{F}_{\lambda_\star}; \cF_\alpha(L)] \leq
C_2(\alpha,
L)\varphi(n),
\]
where the sequences $\lambda(n)$ and $\varphi(n)$ are given in
Table \ref{tab1}, and
constants $C_1$ and $C_2$ are specified in the proof (see
\textup{(A.15)--(A.22)} in \cite{DGJ-suppl}).
%
%
\begin{table}
\tabcolsep=0pt
\caption{The bandwidth order $\lambda(n)$
and the convergence rate of the maximal risk $\varphi(n)$
in the singular and border zones}\label{tab1}
\begin{tabular*}{\tablewidth}{@{\extracolsep{\fill}}lccccc@{}}
\hline
&\multicolumn{3}{c}{\textbf{Border zone} $\bolds{\Theta_{\mathrm b}\mbox{\textbf{:} }\alpha+\beta
=1/2}$} &
\multicolumn{2}{c@{}}{\textbf{Singular zone}
$\bolds{\Theta_{\mathrm s}\mbox{\textbf{:} }\alpha+\beta<1/2}$}\\[-4pt]
& \multicolumn{3}{c}{\hrulefill} & \multicolumn{2}{c}{\hrulefill}\\
&
$\bolds{\beta>1/2}$ & $\bolds{\beta=1/2}$ &$\bolds{\beta<1/2}$
& $\bolds{\alpha+3\beta\geq1/2}$ &
$\bolds{\alpha+3\beta< 1/2}$
\\
\hline
$\lambda(n)$ &
$\frac{n}{\sqrt{\ln n}}$ & $\frac{n}{(\ln n)^{3/2}}$
& $(\frac{n}{\sqrt{\ln n}})^{1/(2\alpha+1)}$
& $n^{{2}/({2\alpha+3-2\beta})}$ &
$n^{{1}/({2\alpha+2\beta+1})}$
\\[6pt]
$\varphi(n)$ & $(\frac{\sqrt{\ln n}}{n})^{\alpha+1/2}$ &
$\frac{(\ln n)^{3/4}}{\sqrt{n}}$ & $\frac{(\ln n)^{1/4}}{\sqrt{n}}$
& $n^{-({2\alpha+1})/({2\alpha+3-2\beta})}$ &
$n^{-({2\alpha+1})/({4\alpha+4\beta+2})}$
\\
\hline
\end{tabular*}
\end{table}
In addition,
if $\lambda_\star=C_3(\alpha, L) \lambda(n/\ln[2\epsilon^{-1}])$, then
for large enough $n$
\[
\Risk_\epsilon[\tilde{F}_{\lambda_\star}; \cF_\alpha(L)] \leq
C_4 (\alpha,
L)\varphi(n/\ln[2\epsilon^{-1}]).
\]
\end{theorem}

Several remarks on the results of Theorems \ref{thupper-bound}--\ref
{thupper-bound-2}
are in order.

\subsubsection*{Remarks}

(1) Theorem \ref{thupper-bound} shows that\vspace*{2pt}
the regular zone $\Theta_{\mathrm r}$ is decomposed into three disjoint regions
with respect to the upper bounds on the risks of $\tilde{F}_{\lambda
_\star}$.
In the zone $\Theta_{{\mathrm r},2}$ where $\beta<1/2$,
the rates of convergence are parametric; because of roughness
of the error density,
here the estimation
problem is essentially a parametric one.
The region $\Theta_{{\mathrm r},1}$ is characterized by nonparametric
rates, while in the border
zone between $\Theta_{{\mathrm r},1}$ and $\Theta_{{\mathrm r},2}$
($\beta
=1/2$) the rate of
convergence differs from the parametric one by a $\ln n$-factor.

(2) The condition on $L$ stated in Theorem \ref{thlower-1} is purely
technical; it requires that the family\vspace*{1pt} $\cF_\alpha(L)$ is
rich enough. It follows from Theorems \ref{thupper-bound} and
\ref{thlower-1} that the estimator $\tilde{F}_{\lambda_\star}$ is
optimal in order in the regular zone if $\alpha>1/2$.

(3) The subdivision of the singular zone $\Theta_{\mathrm s}$ into
two zones
$\Theta_{{\mathrm s},1}=\{(\alpha,\beta)\in\Theta_{\mathrm s}\dvtx
3\beta
+\alpha\geq\frac{1}{2}\}$ and
$\Theta_{{\mathrm s},2}=\{(\alpha,\beta)\in\Theta_{\mathrm s}\dvtx
3\beta
+\alpha<\frac{1}{2}\}$
is a consequence of
two types of upper bounds that we have on the
variance term; see (\ref{eqzone-2}) in \cite{DGJ-suppl}.
In the border zone $\Theta_{\mathrm b}$ the upper bounds on the risk
differ from those
in the regular zone only by logarithmic in $n$ factors.
We do not know if the estimator $\tilde{F}_{\lambda_\star}$ is rate optimal
in the singular and border zones.

(4) Note that the results of
Theorems \ref{thupper-bound} and \ref{thupper-bound-2},
when put together, allow us to establish risk bounds for any pair
$(\alpha,\beta)$ from the
parameter set $\Theta=\{(\alpha,\beta)\dvtx\alpha>-1/2, \beta>0\}$.
In particular,
for any fixed $\alpha>-1/2$, the rate of convergence of the maximal
risk approaches the parametric rate when $\beta$ approaches zero.
We would like to stress the fact that no
tails or moment conditions on $F$ are required to obtain these results;
such conditions
were systematically imposed in the previous work on deconvolution of
distribution functions.

\subsection{Adaptive estimation}\label{secadaptive}

The choice of the smoothing parameter $\lambda$ in
(\ref{eqcdfEstimator1}) is crucial in order to achieve
the optimal estimation accuracy.
As Theorems \ref{thupper-bound}~and~\ref{thlower-1} show, if
parameters\vadjust{\goodbreak}
$\alpha$ and $L$ of the class $\cF_\alpha(L)$ are known, then one
can choose
$\lambda$ in such a way that the resulting estimator is optimal in order.
In practice the functional class $\cF_\alpha(L)$ is hardly known; in these
situations the estimator of Section \ref{secminimax} cannot be implemented.
Note, however, that this does not pose a serious problem in the regular zone
when
$\beta\in(0,1/2)$. Indeed, here if we choose $\lambda=\sqrt{n}$, then
the
resulting estimator will be optimal in order
for any functional class $\cF_\alpha(L)$
satisfying $\lambda_\star=\lambda_\star(\alpha,L) \leq\sqrt{n}$, where
$\lambda_\star$ is defined in Theorem \ref{thupper-bound}.

The situation is completely different in the case $\beta>1/2$.
In this section
we develop an estimator that is nearly optimal for the $\epsilon$-risk
over a
scale of classes
$\cF_\alpha(L)$.
The construction of our adaptive estimator
is based on the general scheme by \cite{lepski}.

\subsubsection{Estimator construction}
Consider
the family of estimators
$\{\tilde{F}_\lambda, \lambda\in\Lambda\}$, where
$\tilde{F}_\lambda$ is defined in (\ref{eqcdfEstimator1}),
$\Lambda:=\{\lambda_j, j=1,\ldots, N\}$ with $\lambda_{\min
}:=\lambda_1$, $\lambda_{\max}:=\lambda_N$,
and
$\lambda_j=2^j\lambda_{\min}$, $j=2, \ldots, N$.
The adaptive\vspace*{1pt} estimator $\tilde{F}$ is obtained
by selection from the family $\{\tilde{F}_\lambda, \lambda\in
\Lambda\}$
according to the following rule.

Let
%
%
\begin{equation}\label{eqc-*}
\omega_1:=\min\{\omega_0, (4b_\zeta)^{-1/\tau}\},\qquad
c_*:= 2\pi^{-2}[2+ (1/\tau)]^2,
\end{equation}
where constants $\omega_0$, $b_\zeta$ and $\tau$
appear in assumption \ref{assE2}.
For any $\lambda\in\Lambda$ we define
%
%
\begin{eqnarray}\label{eqsigma-tilde}
\tilde{\sigma}_\lambda^2 &:=& c_*+
\frac{2}{\pi^2 n}\sum_{j=1}^n \int_{\omega_1}^\lambda\int
_{\omega_1}^\lambda
\frac{1}{\omega\mu}
\Im\biggl\{\frac{e^{i\omega(Y_j-t_0)}}{\hat{f}_\zeta(\omega)}
\biggr\}
\Im\biggl\{\frac{e^{i\mu(Y_j-t_0)}}{\hat{f}_\zeta(\mu)}\biggr\}
\,\rd\omega\,\rd\mu,
\nonumber\\[-8pt]\\[-8pt]
\tilde{\Sigma}_\lambda^2 &:=& \max_{\mu\in\Lambda\dvtx\mu\leq
\lambda}
\tilde{\sigma}_\mu^2 .
\nonumber
\end{eqnarray}
Note that $\tilde{\sigma}_\lambda^2$ can be computed from the data
(the parameters $\tau$ and $\omega_1$ are determined completely by
$\hat{f}_\zeta$; hence they are known).
In fact, $\tilde{\sigma}^2_\lambda n^{-1}$ is a plug-in
estimator of an upper bound on the variance of $\tilde{F}_\lambda$, while
$\tilde{\Sigma}_\lambda^2$ is a ``monotonization'' of
$\tilde{\sigma}^2_\lambda$
with respect to $\lambda$.

Define
\[
\tilde{v}_\lambda^2 :=
\tilde{\Sigma}_\lambda^2 +
11\bar{m}^2 \lambda^{2\beta}n^{-1} \ln(4N^2\epsilon^{-1}),\qquad
\lambda\in\Lambda,
\]
where
%
%
\begin{equation}\label{eqm-bar}
\bar{m}:=
\sqrt{2c_*} + (\pi c_\zeta\beta)^{-1} 2^{1+(\beta/2-1)_+} [2+
\beta\ln_+
(1/\omega_1)],
\end{equation}
and constant $c_\zeta$ appears in assumption \ref{assE1}.

Let $\vartheta:=\sqrt{2} (\sqrt{2}-1)^{-1}[1+
\sqrt{3\ln(4N\epsilon^{-1})}]$; then
with every estimator
$\tilde{F}_\lambda$, $\lambda\in\Lambda$ we associate the interval
%
%
\begin{equation}\label{eqQ}
Q_\lambda:= [\tilde{F}_\lambda- \vartheta\tilde{v}_\lambda
n^{-1/2},
\tilde{F}_\lambda+ \vartheta\tilde{v}_\lambda n^{-1/2}].
\end{equation}
Define
%
%
\begin{equation}\label{eqlambda-choice}
\tilde{\lambda}:=\min\biggl\{\lambda\in\Lambda\dvtx\bigcap_{
\mu\geq\lambda, \mu\in\Lambda} Q_\mu\ne\varnothing\biggr\},
\end{equation}
and set finally
%
%
\begin{equation}\label{eqest}
\tilde{F}:=\tilde{F}_{\tilde{\lambda}}.
\end{equation}
Note that $\tilde{\lambda}$ is well defined: the intersection
in (\ref{eqlambda-choice}) is nonempty for $\lambda=\lambda
_{\max}$.
%

\subsubsection{Oracle inequality}
We will show that the estimator $\tilde{F}$ mimics the oracle estimator
$\tilde{F}_{o}$ which is defined as follows:

Let
\begin{eqnarray*}
\sigma_\lambda^2 &:=& c_* + \frac{2}{\pi^2}
\bE\biggl[\int_{\omega_1}^\lambda\frac{1}{\omega}
\Im\biggl\{\frac{e^{i\omega(Y_j-t_0)}}{\hat{f}_\zeta(\omega)}
\biggr\} \,\rd\omega
\biggr]^2,
\\
\Sigma_{\lambda}^2 &:=& \max_{\mu\in\Lambda\dvtx\mu\leq\lambda}
\sigma_\mu^2,\qquad
\lambda\in\Lambda.
\end{eqnarray*}
It is shown in the proof of
Lemma \ref{lemvariance} (see Section A.1.2 in \cite
{DGJ-suppl})
that $\sigma_\lambda^2n^{-1}$ is an upper bound on the variance
of the estimator $\tilde{F}_\lambda$ associated with parameter
$\lambda$. Note that $\tilde{\sigma}^2_\lambda$
defined in (\ref{eqsigma-tilde}) is the empirical counterpart
of the quantity $\sigma_\lambda^2$.
Define
\[
v_\lambda^2 :=\Sigma_\lambda^2 + 11\bar{m}^2\lambda^{2\beta}
n^{-1} \ln(4N^2\epsilon^{-1}).
\]
Given $\alpha>0$ and $L>0$ let
\[
\lambda_o=\lambda_o(\alpha, L):= \min\bigl\{ \lambda
\in\Lambda\dvtx v_\lambda n^{-1/2} \geq
2\sqrt{2}\pi^{-1/2} L \lambda^{-\alpha-1/2}
\bigr\}
\]
and define $\tilde{F}_o:=\tilde{F}_{\lambda_o}$.

The oracle estimator $\tilde{F}_o$ has attractive minimax properties
over classes $\cF_\alpha(L)$.
In particular, it is easily verified that for any class $\cF_\alpha
(L)$ such that
$\lambda_o\leq[11\bar{m}^2 \ln(4N\epsilon^{-1})]^{-1}n$
one has
\[
\Risk_{\epsilon}[\tilde{F}_o; \cF_\alpha(L)] \leq
2v_{\lambda_o} n^{-1/2} \leq\kappa_1 \psi_{n, \epsilon}(\alpha
, L) +
\kappa_2\phi_{n,\epsilon}.
\]
Here $\psi_{n,\epsilon}$ is the upper bound of Theorem \ref{thupper-bound}
on the risk of the estimator
$\tilde{F}_{\lambda_\star}$ that ``knows'' $\alpha$ and $L$, $\phi
_{n,\epsilon}$ is defined
in Theorem \ref{thlower-1}, and $\kappa_1$ and $\kappa_2$ are
constants independent of
$\alpha$ and $L$.
Thus, the risk of the oracle estimator admits the same upper bound as the
risk of the estimator
$\tilde{F}_{\lambda_{\star}}$ that is based on the knowledge of the
class parameters $\alpha$ and $L$.

Now we are in a position to state a bound on
the risk of the estimator~$\tilde{F}_{\tilde{\lambda}}$.

\begin{theorem}\label{thadaptive}
Suppose that assumptions \ref{assE1}, \ref{assE2} hold, $\beta
>1/2$ and let
\[
\lambda_{\max}=[11\bar{m}^2 \ln(4N\epsilon^{-1})]^{-1} n.
\]
If $\tilde{F}_{\tilde{\lambda}}$ is the estimator defined
in (\ref{eqQ})--(\ref{eqest})
then for any class $\cF_\alpha(L)$ with $\alpha>0$
such that $\lambda_{\min}\leq\lambda_o(\alpha, L)\leq\lambda
_{\max}$,
one has
\[
\Risk_{\epsilon} [\tilde{F}_{\tilde{\lambda}}; \cF_\alpha(L)]
\leq\bigl(3- 1/\sqrt{2}\bigr)
\vartheta v_{\lambda_o}n^{-1/2}.
\]
\end{theorem}

Estimator (\ref{eqQ})--(\ref{eqest}) attains the optimal rates of
convergence
with respect to $\epsilon$-risk
within a $\ln(N\epsilon^{-1})$-factor\vadjust{\goodbreak}
over the collection of functional
classes $\cF_\alpha(L)$. In particular,
if $\lambda_{\min}$ is chosen to be a constant, and $\lambda_{\max
}\asymp n^l$
for some $l\geq1$, then
$N=\ln(\lambda_{\max}/\lambda_{\min})/\ln2 \asymp\ln n$, and
the $\epsilon$-risk of the adaptive estimator $\tilde{F}_{\tilde
{\lambda}}$
is within a $\ln\ln n$-factor of the minimax $\epsilon$-risk
for a scale of Sobolev classes. It can be shown that
this $\ln\ln n$-factor is unavoidable price for adaptation
when the accuracy is measured by the $\epsilon$-risk; see, for
example, \cite{spokoiny}.

%
\section{Minimax and adaptive affine estimation in discrete
deconvolution model}
\label{secminimax-linear}

The results of Section \ref{secminimax} imply that in the
regular zone the minimax rates of
convergence on the Sobolev classes are attained by linear estimator
(\ref{eqcdfEstimator1}). It seems interesting to compare the
performance of
estimator (\ref{eqcdfEstimator1}) and its adaptive version in
Section \ref{secadaptive} with that of the minimax linear estimator.
%

Consider the estimation problem as follows; cf.
\cite{JudNem}, Problem 2.2:
\renewcommand{\theproblem}{D}
\begin{problem}\label{probD}
We observe $n$ independent realizations
$\eta_1,\ldots,\eta_n$ of a random variable $\eta$, taking values in
$\bS=\{1,\ldots,m\}$.
The distribution of $\eta$ is identified with a vector $p$
from the
$m$-dimensional simplex $\P_{m}=\{y\in\bR^{m}\dvtx y\ge0,\sum_iy_i=1\}
$ by setting
$p_k=\bP\{\eta=k\}$, $1\leq k\leq m$.
Suppose that vector $p$
is affinely parameterized
by an $M$-dimensional
``signal''-vector of unknown ``parameters'' $x\in
\mathcal{X}\subset\P_M\dvtx p=Ax=[[Ax]_1;\ldots;[Ax]_m]$. Here $Ax$ is
the linear
mapping with $A\mathcal{X}\subset\P_{m}$, and $[a]_j$ stands for the $j$th
element of $a$.
Our goal is to estimate a given linear form $g(x)=g^Tx$
at the point $x$ underlying the observation $\eta^n$.
\end{problem}

It is obvious that if distributions of $X$ and $\zeta$ are compactly
supported, or can be ``compactified'' (i.e., for any $\eep>0$ one can
point out
bounded intervals of probability $1-\eep$ for $X$ and $\zeta$), then
under very minor regularity conditions on $f_\zeta$ and $F$, the
Problem \ref{probD}
approximates the initial distribution deconvolution problem with ``arbitrary
accuracy.''
The latter means that given $\eep>0$ we can compile the discretized
problem such that its $\delta$-solution is the solution to the initial
continuous problem with the accuracy $\delta+\eep$
with probability $1-\eep$.

We consider
the following discretization of the deconvolution problem:
\begin{longlist}[(4)]
\item[(1)] Let $ J=[a_0,a_m]$
be the (finite) observation domain, and let $a_0<a_1<a_2<\cdots<a_{m-1}<a_m$.
We split $ J$ into $m$ intervals
$ J_1=[a_0,a_1], J_2=(a_1,a_2],\ldots, J_m=(a_{m-1},a_m]$. We denote
$p_k=\rP\{Y\in J_k\}$, $k=1,\ldots,m$.
\item[(2)] Suppose that the (finite) interval $I=[b_0,b_M]$ contains
the support of all $F\in\cF$. Let $b_0<b_1<b_2<\cdots<b_{M}$, we
partition $I$ into $M$ intervals
$I_1=[{b}_0,b_1], I_2=(b_1,b_2],\ldots,I_M=(b_{M-1},b_M]$.
We denote $x_k=\rP\{X\in I_k\}$, $k=1,\ldots,M$.
\item[(3)] Denote $\bar{b}_k=(b_{k-1}+b_k)/ 2$. Define the $m\times
M$ matrix $A=(A_{jk})$ with elements
\begin{eqnarray*}
A_{jk}&=&\bP\{{\bar{b}_k}+\zeta\in J_j\}
\\
&=&
\cases{
\bP\{a_{0}-\bar{b}_k\le\zeta\le a_{1}-\bar{b}_k\}
,&\quad$k=1,\ldots,M, j=1$,\vspace*{2pt}\cr
\bP\{a_{j-1}-\bar{b}_k<\zeta\le a_{j}-\bar{b}_k\},
&\quad$k=1,\ldots,M, j=2,\ldots,m$,}
\end{eqnarray*}
and the vector $g=g(t_0)\in\bR^M$, with $g_k=\ind(\bar{b}_k\le t_0)$,
$k=1,\ldots, M$.
The elements $A_{jk}$ of $A$ are the approximations of conditional probabilities
$\rP\{Y\in J_j|X\in I_k\}$, and
$g^Tx$ is an approximation of $F(t_0)$.
\item[(4)] Consider discrete observations $\eta_i\in\{1,\ldots,m\}$ as follows:
\[
\eta_i=\ind(a_{0}\le Y_i\le a_1)+\sum_{j=2}^m j \cdot\ind
(a_{j-1}< Y_i\le a_j),\qquad
i=1,\ldots,n.
\]
\end{longlist}
%
If the sets $I$ and $J$ are selected so that
$\bP\{X\in I\}\ge
1-\eep$, $\bP\{Y\in J\}\ge1-\eep$ for any $F\in\cF$, if
$\cF$ is
the class of ``regular distributions'' and the noise distribution 
possesses some regularity, and if the partitions of $I$ and $J$ are
``fine enough,''
then solving Problem \ref{probD} with
$\mathcal X$ being the corresponding $M$-dimensional
cross-section of $\cF$ will provide us with an estimation $\tilde{g}$ of
$F(t_0)$ in the continuous deconvolution problem.

We now
concentrate on solving the deconvolution problem
in the discrete model.

\subsection{Minimax estimation in the discrete model}
\label{sectdiscrmmax}

An estimate of $g(x)$---a candidate solution to our problem---is a
measurable function
$\wt{g}=\widetilde{g}(\eta^n)\dvtx\bS^n\to\bR$. Given tolerance
$\e\in(0,1)$, we define the \textit{$\e$-risk} of such an
estimate on $\mathcal{X}$ as
\[
\Risk_\e(\wt{g};\mathcal{X})=\inf\Bigl\{\delta\dvtx\sup_{x\in
\mathcal{X}}\bP_{x}\{|\wt{g}(\eta^n)-g^Tx|>\delta\}<\e
\Bigr\},
\]
where
$\bP_{x}$ stands for the distribution of
observations $\eta^n$ associated with the ``signal''~$x$.
The minimax optimal $\e$-risk is
\[
\Risk^*_{\e}(\mathcal{X})=\inf_{\wt{g}(\cdot)}
\Risk_{\e}(\wt{g};\mathcal{X}).
\]

We are particularly interested in the family of estimators of the
following structure:
\[
\wt{g}_{\vf,c}(\eta^n)={1\over n}\sum_{i=1}^n \vf(\eta_i)+c
={1\over n}
\sum_{i=1}^n \sum_{k=1}^m \vf_k \ind(\eta_i=k) +c.
\]
We refer to such estimators $\wt{g}_\vf$ as \textit{affine}.
In other words, $\wt{g}_\vf$ is
an affine function of empirical distribution: for some $\vf\in\bR^m$
and $c\in\bR$,
\[
\wt{g}_{\vf,c}(\eta^n)=\sum_{k=1}^m \vf_k \wt{P}_n(k)+c,
\]
where
$\wt{P}_n$ is the empirical distribution of the observation sample
$
\wt{P}_n(k)={1\over n}\sum_{i=1}^n \ind(\eta_i=k).
$
An important property of the class of affine estimators, when applied
to Problem \ref{probD} with convex set $\mathcal{X}$, is that one can choose
an estimator from the class such that its $\e$-risk attains (up to a
moderate $\mbox{constant}\approx2$; see Theorem~\ref{thethemain1} below)
the minimax $\e$-risk $\Risk_\e^*(\cX)$.

From now on let us assume that $\mathcal{X}\subset\bR^M$ is a convex
closed (and, being a subset of an $M$-dimensional simplex, compact) set.

Let us consider the affine estimator $\wt{g}_\e$ of $g^Tx$
\[
\wt{g}_\e(\eta^n)\equiv
\wt{g}_{\bar{\varphi},\bar{c}}(\eta^n)=\sum_{k=1}^{m}\bar{\vf
}_{k}\wt{P}
_n(k)+\bar{ c },
\]
%
in which the parameters $\bar{\vf}$ and $\bar c$ of $\wt{g}_\e$ are
defined as
follows.

Consider the optimization problem
%
%
\begin{eqnarray}\label{eqcon2}
\overline{S}(\e)&=&\max_{x,y\in\mathcal{X}}
\Biggl\{\half g^T(y-x),\nonumber\\[-8pt]\\[-8pt]
&&\hphantom{\max_{x,y\in\mathcal{X}}
\Biggl\{}h(x,y;\e)\equiv n\ln\Biggl(\sum_{j=1}^{m}\sqrt{[Ax]_j[Ay]_j}
\Biggr)+\ln(2/\e)\ge0\Biggr\}.\nonumber
\end{eqnarray}
%
Let $(\bar{x},\bar{y})$ be an optimal solution to (\ref{eqcon2}),
and let
$\nu\ge0$
be the Lagrange
multiplier of the constraint $h(x,y;\e)\ge0$. We set
\[
\bar{c}=\half g^T[\bar{y}+\bar{x}] 
\quad\mbox{and}\quad
\bar\vf_{j}=\nu n\ln\Biggl[\sqrt{[A\bar{y}]_j\over[A\bar
{x}]_j}\Biggr],\qquad
j=1,\ldots,m.
\]

We have the following result.
\begin{theorem}\label{thethemain1}
Let $\e\in(0,1/4]$. Then the $\epsilon$-risk of the estimator
$\wt{g}_{\epsilon}$ satisfies
%
%
\begin{equation}\label{eqteq0}
\Risk_\e(\wt{g}_\epsilon;\mathcal{X})\le\overline{S}(\epsilon
)\le
\vartheta(\epsilon)\Risk^*_\e(\mathcal{X}),\qquad
\vartheta(\epsilon)={2\ln(2/\epsilon)\over
\ln[1/(4\epsilon)]} .
\end{equation}
\end{theorem}

Note that $\vartheta(\epsilon)\to2$ as $\epsilon\to0$; thus for small
tolerance levels the $\epsilon$-risk of the estimator $\tilde
{g}_\epsilon$ is within factor $\approx2$ of the minimax
$\epsilon$-risk. It is important to emphasize that
$\widetilde{g}_{\epsilon}$ is readily given by a solution to the
explicit convex program (\ref{eqcon2}), and as such, it can be found in
a computationally efficient fashion, provided that $\mathcal{X}$ is
computationally tractable.

In the ``historical perspective'' the affine estimator $\wt{g}_\e$
represents an alternative to the binary search estimator $\wt{g}_B$,
proposed in \cite{DonLiuII} for the case of ``direct'' observations. It
can be shown that the $\e$-risk $\Risk_\e(\wt{g}_B;\mathcal {X})$ of
that estimator satisfies $\Risk_\e(\wt{g}_B;\mathcal{X})\le C \Risk
^*_\e(\mathcal{X})$ for small $\e$ (e.g., one can prove that
$C\le26$ whenever $\epsilon\le 0.01$). To the best of our knowledge,
risk bound (\ref{eqteq0}) in Theorem~\ref{thethemain1} for the
estimator $\wt{g}_\e$ is much better than those available for the
binary search estimator.

Note that the constraint $h(x,y;\epsilon)\ge0$ of the problem (\ref
{eqcon2}) can be rewritten as follows:
\[
\rho(x,y)\ge(\epsilon/2)^{1/n},
\]
where
\[
\rho(x,y)= \sum_{k=1}^m \sqrt{[Ax]_k[Ay]_k}
\]
is the \textit{Hellinger affinity} of distributions $A(x)$ and $A(y)$;
cf. \cite{LeCam73} and \cite{LeCam86}, Chapter~4.
Thus the optimal value $\overline{S}(\epsilon)$ of the optimization
problem (\ref{eqcon2}) can be seen as modulus of continuity of the linear
functional $g(\cdot)$ over the class $\mathcal{X}$ of distributions
``with respect to Hellinger affinity.'' If
$\on\ln[1/\epsilon]=o(1)$ we have $\rho(x,y)\approx1$ and
\[
H^2(x,y)=1-\rho(x,y)\approx-\ln\rho(x,y),
\]
where $H(x,y) $ is the \textit{Hellinger distance} between $x$ and $y$. In
this limit we have
\[
\overline{S}(\epsilon)\approx\half\omega\Biggl(\sqrt{\ln
[2/\epsilon]\over
2N}\Biggr)\equiv
\max_{x,y\in\mathcal{X}}\Biggl\{{\half} g^T(y-x), H(x,y)\le
\sqrt{\ln[2
/\epsilon]\over2n}\Biggr\}.
\]
%
Here $\omega(\cdot)$ is the ``modulus of continuity of $g$ over
$\mathcal{X}$
with respect to Hellinger distance,'' introduced in \cite{DonLiuII}.
Therefore, bound (\ref{eqteq0}) can be seen as a finite-dimensional
nonasymptotic counterpart
of \cite{DonLiuII}, Theorem 3.1.

\subsection{Adaptive version of the estimate}\label{secadapt2}
\label{sectadaptive}

Consider a modification of our estimation problem where the set
$\mathcal{X}$, instead of
being given in advance, is known to be one of the sets from the
collection of nonempty convex compact sets $\mathcal{X}^1, \mathcal
{X}^2,\ldots,
\mathcal{X}^N$ in $\bR^M$. We aim to construct an \textit{adaptive
estimator} of the
linear form $g^Tx$, given that $x$ is an element of some $\mathcal
{X}_i$ in the
collection. Here we consider the simple case where the sets are nested.
$\mathcal{X}^1\subset\mathcal{X}^2\subset\cdots\subset\mathcal
{X}^N$. Note that in the case of
nonnested sets
an adaptive estimator can be constructed following the ideas of \cite
{CaiLow2004}.

Given a linear form $g^Tz$ on $\bR^M$, let
$\Risk^k(\wt{g})$ and $\Risk^k_*$ be,
respectively, the $\epsilon$-risk of an estimate $\wt{g}$ on
$\mathcal{X}^k$, and the minimax optimal $\epsilon$-risk of recovering
$g^Tx$ on $\mathcal{X}^k$. Let also $S_k(\cdot)$ be the function
$\overline{S}(\cdot)$ in
(\ref{eqcon2}) associated with $\mathcal{X}=\mathcal{X}^k$. As it is
immediately seen,
the functions $S_k(\cdot)$ grow with $k$. Our goal is to modify
the estimate $\wt{g}$ we have built in such a way that the
$\epsilon$-risk of the modified estimate on $\mathcal{X}^k$ will be
``nearly'' $\Risk^k_*$ \textit{for every} $k\leq N$. This
goal can be achieved by a straightforward application of
Lepski's adaptation scheme as follows.

Given $\epsilon>0$,
let
$\wt{g}^k(\cdot)$ be the affine estimate with the
$(\epsilon/N)$-risk on $\mathcal{X}^k$ not exceeding
$S_k(\epsilon/N)$ as provided by Theorem \ref{thethemain1} which is
applied with
$\epsilon/N$ substituted for $\epsilon$
and $\mathcal{X}^k$ substituted for $\mathcal{X}$. Then
\[
\sup_{x\in
\mathcal{X}^k}\bP_x\{|\wt{g}^k(\eta^n)-g^Tx|>S_k(\epsilon
/N)\}\leq
\epsilon/N\qquad \forall k\leq N.
\]
Given observation $\eta^n$, let us say that the index $k\leq N$ is
\textit{$\eta^n$-good}, if for all $k'$ satisfying $k\leq k'\leq N$ one has
\[
|\wt{g}^{k'}(\eta^n)-\wt{g}^k(\eta^n)|\leq
S_k(\epsilon/N)+S_{k'}(\epsilon/N).
\]
Note that $\eta^n$-good indices do exist (e.g., $k=N$). Given $\eta^n$,
we can find the smallest $\eta^n$-good index $k=k(\eta^n)$; our estimate
is nothing but $\wt{g}(\eta^n)=\wt{g}^{k(\eta^n)}(\eta^n)$.
\begin{proposition}\label{adaptiven}
Assume that $\epsilon\in(0,1/4)$, and let
\[
\vartheta=3{\ln(2N/\epsilon)\over\ln(2/\epsilon)}.
\]
Then
\[
\sup_{x\in\mathcal{X}^k}
\bP_x\{|\wt{g}(\eta^n)-g^Tx|>\vartheta
S_k(\epsilon)\}<\epsilon\qquad \forall(k, 1\leq k\leq N);
\]
whence also
\[
\Risk^k(\wt{g})\leq
{6\ln(2N/\epsilon)\over\ln[1/
(4\epsilon)]}\Risk^k_*\qquad \forall(k,1\leq k\leq N).
\]
\end{proposition}

The proof of the proposition follows exactly same steps as that of
Proposition~5.1 of \cite{JudNem}, and it is omitted.

\section{Numerical examples}\label{secnumeric}

To illustrate our results we present here examples of implementation of the
adaptive estimation procedures of
Sections \ref{secadaptive} and \ref{sectadaptive}.

We consider three measurement error distributions scenarios:
\begin{longlist}[(iii)]
\item[(i)] Gamma distribution $\Gamma(0,2,1/(2\sqrt{2}))$
with the shape parameter $2$ and the scale $\frac{1}{2\sqrt{2}}$
(the standard
deviation of the error is equal to $0.5$). Here $\Gamma(\mu,\alpha
,\theta)$
stands for the Gamma distribution with location $\mu$, shape
parameter $\alpha$ and scale $\theta$, such that its density is
$
[\Gamma(\alpha)\theta^\alpha]^{-1} (x-\mu)^{\alpha-1}
\exp\{-(x-\mu)/\theta\}\ind(x\ge\mu)$.\vspace*{1pt}
\item[(ii)] Mixture\vspace*{1pt} of Laplace distributions
$\half\mathcal{L}(-1,\half)+\half\mathcal{L}(1,\half)$; here
$\mathcal{L}(\mu, a)$
stands for the Laplace distribution with the density
${(2a)^{-1}}e^{-{|x-\mu|/a}}$.
\item[(iii)] Normal mixture $\half\mathcal{N}(0,{\small1\over
4})+\half\mathcal{N}(2,{\small1\over4})$.
\end{longlist}
We
consider
three distributions of $X$:
\begin{longlist}[(3)]
\item[(1)] mixture of ``shifted'' Gamma distributions:
$0.3 \Gamma(0,0.5,2)+0.7 \Gamma(5,0.5,2)$;
\item[(2)] mixture of Laplace distributions
$0.3 \mathcal{L}(-1.5,0.5)+0.7 \mathcal{L}(1.7,0.25)$;
\item[(3)] normal mixture
$0.6 \mathcal{N}(0.15827,1)+0.4 \mathcal{N}(1,0.0150)$.
\end{longlist}

Note that in the case (i) of $\Gamma(0,2,\theta)$ error distribution
the estimator (\ref{eqcdfEstimator1}) can be
computed explicitly: we have
$
\tilde{F}_\lambda={1\over2}-{1\over\pi n}\sum_{i=1}^n I_\lambda(Y_i-t_0),
$
where
\[
I_\lambda(y)=\Si(\lambda y) +
y^{-1}[\theta^2\lambda\cos(\lambda y)-2\theta\sin(\lambda y)]
- y^{-2}\theta^2\sin(\lambda y),
\]
and $\Si(x)=\int_0^{x} \omega^{-1} \sin\omega\,\rd\omega$ is the sine
integral
function.
Then the adaptive estimation algorithm of Section \ref{secadaptive}
is implemented for the grid $\Lambda=\{\lambda\in[0.01\dvtx 0.05\dvtx 10]\}$.

Estimation procedures, described in Section \ref{sectdiscrmmax}, were
implemented using \texttt{Mosek} optimization software \cite{mosek}. The
observation space and the signal space were split into $m=M=200$ bins.
{The adaptation procedure was implemented over 17 linear
estimators corresponding to the classes
$\mathcal{X}^1,\ldots,\mathcal{X}^{17}$ of ``Lipschitz-continuous''
discrete distributions with Lipschitz constants on the geometric grid, scaled
from 0.001 to 1 [if reduced to continuous densities, it corresponds to
the approximate range
of Lipschitz constant from $O(0.1)$ to $O(100)$].}

\begin{figure}

\includegraphics{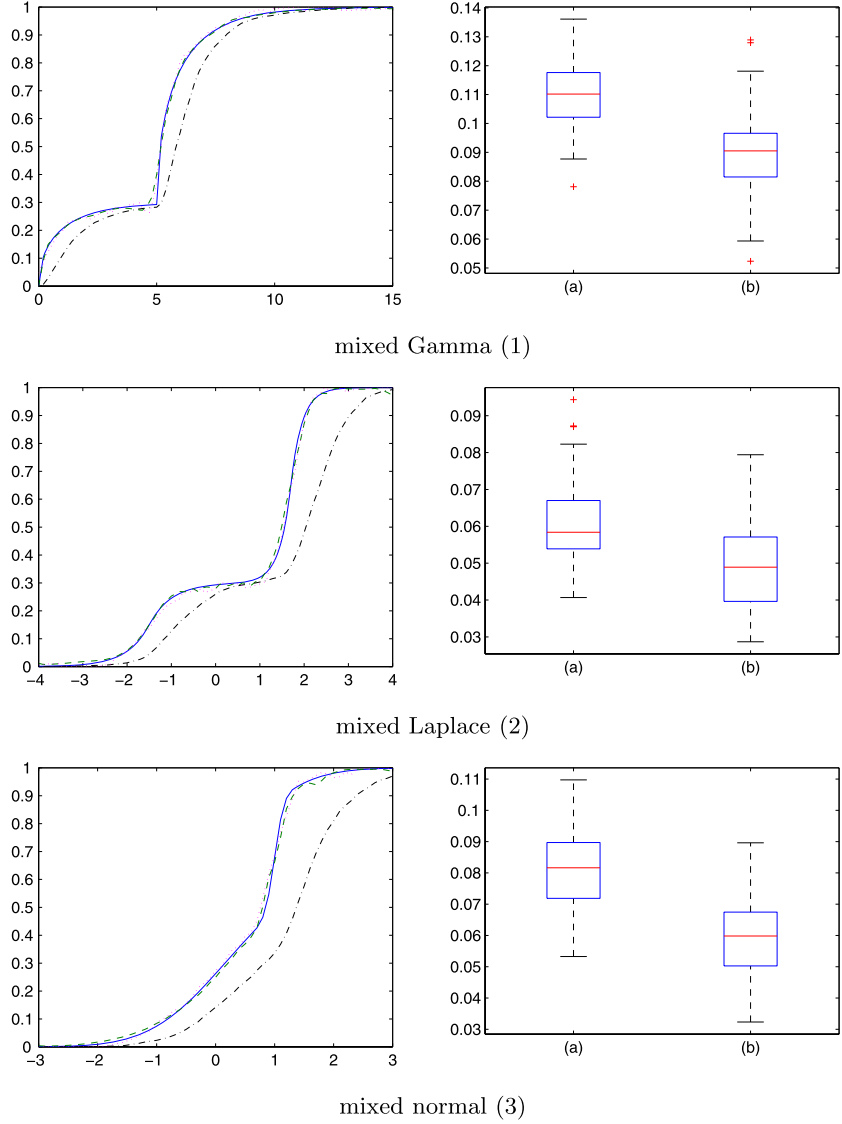}

\caption{Simulation results for the Gamma error scenario. On the left:
true cdf (solid line),
adaptive estimator $\wt{g}(\eta^n)$ of Section \protect\ref{sectadaptive} (dashed line), adaptive
estimator $\wt{F}_{\wt{\lambda}}$ of Section \protect\ref{secadaptive} (dotted line)
and
the edf
of the observations (dash--dot line).
On the right:
the boxplots of the maximal estimation error of $\wt{g}(\eta^n)$
\textup{(a)}
and $\wt{F}_{\wt{\lambda}}$ \textup{(b)}.}
\label{fig2}
\vspace*{-12pt}
\end{figure}

\begin{figure}

\includegraphics{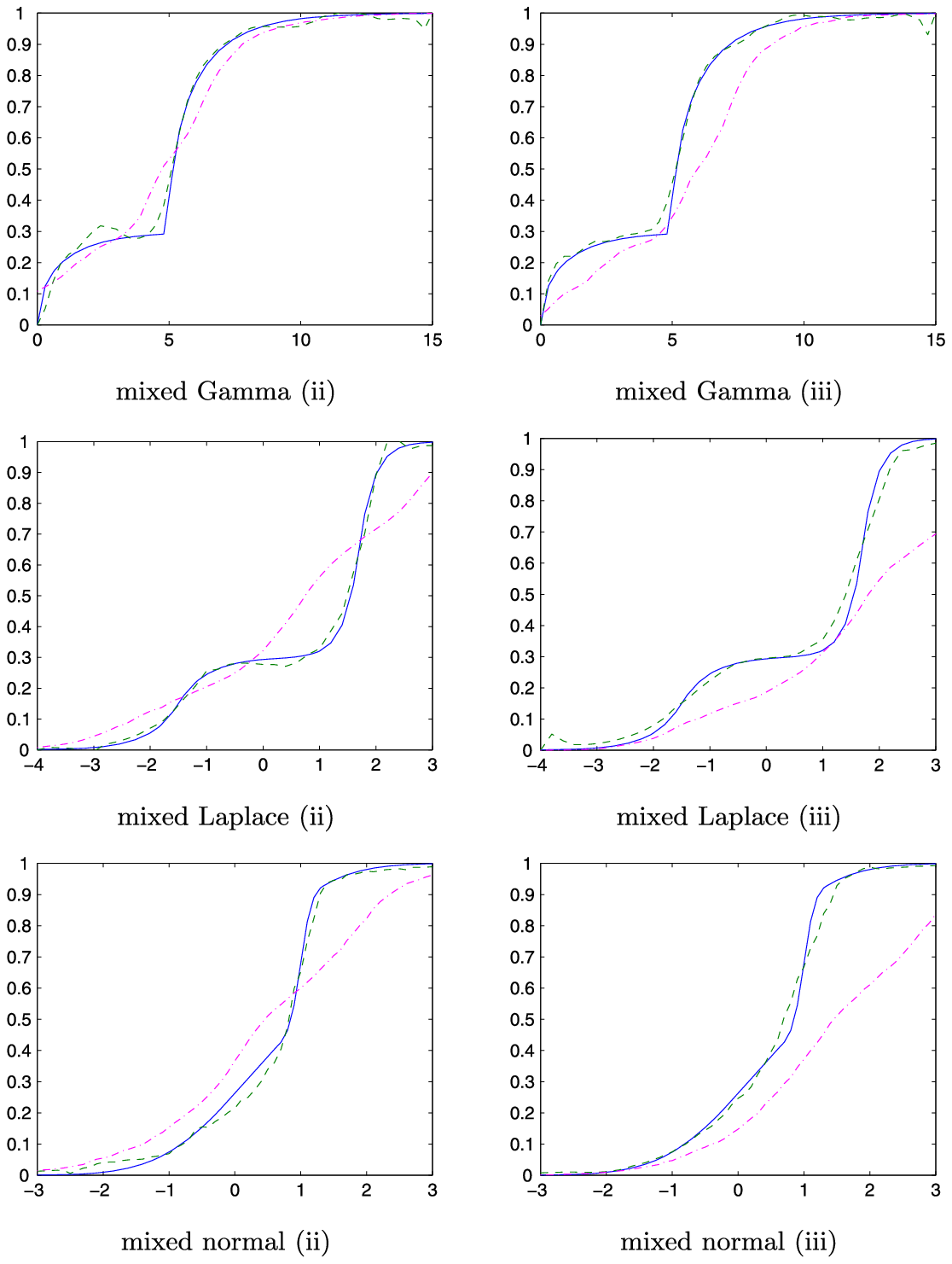}

\caption{Simulation results: true cdf (solid line), adaptive
estimator (dashed line) and empirical distribution function of the observation
(dash--dot line).
On the left, \textup{(ii)} are the results for mixed Laplace noise; on
the right, \textup{(iii)} are the results for the mixed normal noise.} \label{fig10}
\vspace*{-12pt}
\end{figure}

\begin{figure}

\includegraphics{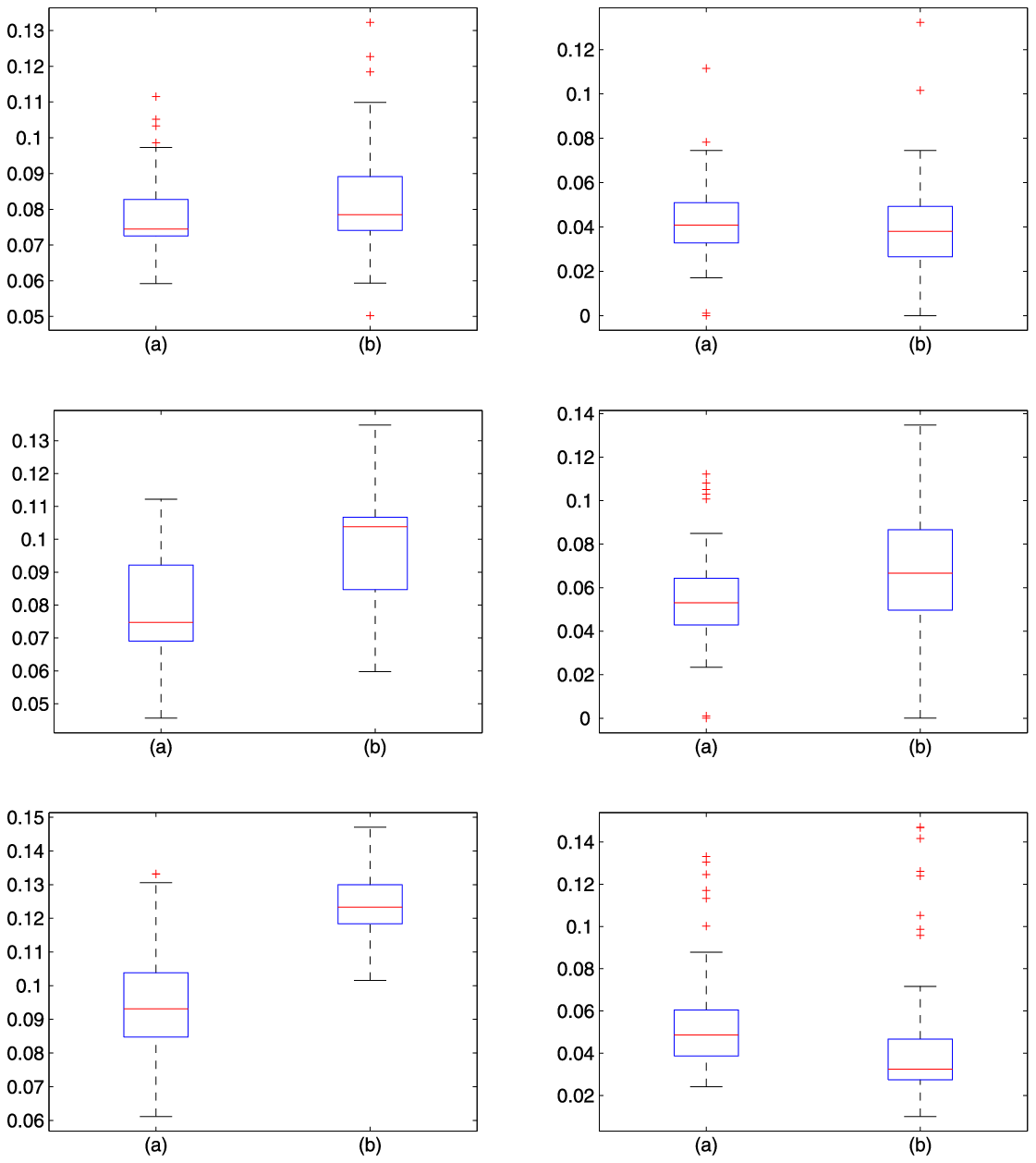}

\caption{Estimation error distribution. Left column: empirical
distribution of the maximal error of estimation over a regular grid;
right column: distribution of the maximal over 100 simulations
estimation error over the points of the grid. On each plot the left
boxplot \textup{(a)} corresponds to the mixed Laplace noise, while the
right boxplot \textup{(b)} corresponds to the mixed normal noise.} \label{fig20}
\vspace*{-12pt}
\end{figure}

The simulation has been repeated for $100$ observation samples
of size $n=2\mbox{,}000$.
On Figure \ref{fig2} we present simulation results for the scenario
%
%
(i) when the error distribution follows the $\Gamma(0,2,1/(2\sqrt{2}))$
law. The left column displays ``typical'' results of estimation
corresponding to three signal distributions. We present the true
distribution (solid line), the estimate $\wt{F}_{\wt{\lambda}}$ of
Section \ref{secadaptive} (dotted line), the estimate $\wt{g}(\eta^n)$
of Section \ref{sectadaptive} (dashed line) and the empirical
distribution of the observations (dash--dot line). The boxplots on the
right display resume the corresponding empirical distributions of the
maximal estimation error over 50 points of the regular grid on the
support of $f$ for two estimators: (a) for $\wt{g}(\eta^n)$ of
Section~\ref{sectadaptive} and (b) for the $\wt{F}_{\wt{\lambda}}$ of
Section \ref{secadaptive}. On Figure \ref{fig10} we present ``typical''
results for adaptive estimator $\wt{g}(\eta^n)$ of Section
\ref{sectadaptive} under the error scenarios (ii) (on the left) and
(iii) (on the right). Similarly to Figure \ref{fig2} we plot true cdf
(solid line), adaptive estimator $\wt{g}(\eta^n)$ of Section
\ref{sectadaptive} (dashed line) and the observation edf (dash--dot
line).
%
%
The results of this simulation are summarized on Figure \ref{fig20}.
The first
boxplot (the left column plots) represents the distribution of the maximal
estimation error over 50 points of the regular grid on the support of $f$.
%
%
Next, for each point in the grid we compute the maximal estimation
error over 100 simulations, the distribution of maximal errors ``over
the points of the grid'' is represented on the second boxplot (plots on
the right column).

\subsection*{Remarks}
The numerical examples in this section illustrate strong and weak
points of the proposed estimators related to practical implementation.
They can be summarized as follows.

The adaptive estimator of Section \ref{secadaptive}
is based on the choice of the unique smoothing parameter $\lambda$.
This imposes a ``natural'' family of nested classes
and facilitates implementation of the adaptation scheme.
Yet, this estimator should be ``explicitly
tuned'' for a specific distribution of the errors.
In particular,
the integral computation in (\ref{eqcdfEstimator1})
for a given distribution of $\zeta$
may become very tedious.
Even though our theoretical
results are proved
under the
condition that $|\hat{f}_\zeta(\omega)|\ne0$ for all
$\omega\in\bR$, in practical implementation
the estimator (\ref{eqcdfEstimator1}) could
be modified in order to allow characteristic functions $\hat{f}_\zeta$
vanishing
at finite number of points
on $\bR$. In this case
the integration domain in (\ref{eqcdfEstimator1})
should exclude some properly specified
vicinities of the points where $\hat{f}_\zeta$ vanishes.

In contrast to this, the
adaptive estimator in Section \ref{secadapt2}
can be easily tuned to \textit{any
noise distribution} and convex target distribution class.
For instance,
the characteristic function of noise in the Laplace scenario (ii)
vanishes at
some points, what
precludes the possibility of utilizing
the estimator of Section \ref{secadaptive} without proper modifications.
Note that one can easily incorporate any additional available
information on the unknown distribution that
can be expressed as a convex constraint in the corresponding optimization
problem. The typical examples of such constraints are
unimodality, symmetry, monotonicity and moment bounds.
However, this freedom comes at a price:
the family $\cX^1\subset
\cdots\subset\cX^N$ of
the embedded classes for the
adaptive estimator in Section \ref{secadapt2} should be constructed ``by
hand.''
The computation of
the adaptive affine estimator of
Section \ref{sectadaptive} is also a
heavy numerical task. In particular, in our setting it involves solving $17$
conic quadratic optimization problems with 1,006
variables, $809$ linear and $202$ conic constraints.

It is well known that the normal noise in the deconvolution problem
results in a
very poor quality of estimation \cite{Fan}.
In particular, the minimax rate of convergence in this case is
$O((\ln n)^{-\gamma})$ with $\gamma>0$ depending on
the exponent $\alpha$ of the regularity class $\mathcal{F}_\alpha(L)$.
Fortunately, these
pessimistic results are concerned with
\textit{the asymptotic} as $n\to\infty$ behavior of
the estimators. We observed that the estimation procedures exhibit much better
performance for
small or moderately sized observation samples. On the other hand, this
performance does not improve when the sample size grows up: in our
experiments, for instance, the estimation accuracy, measured by
$\ell_\infty$-error over a regular grid in the distribution domain, improved
only by the factor $\approx2$ when we increased the sample size from
$n=2\mbox{,}000$ to $n=100\mbox{,}000$.\

\section{Proofs}
\label{secproofs}

This section is organized as follows.
In Section \ref{secproof-1-2} we state main results
that are used in the proof of Theorems \ref{thupper-bound}
and \ref{thupper-bound-2} and briefly discuss the proof outline.
Then in Section \ref{subsecadaptive00} we prove Theorem \ref{thadaptive}.
Full proofs of all auxiliary results and additional technical details
are given in the supplementary
paper~\cite{DGJ-suppl}.

%
\subsection{\texorpdfstring{Proofs of Theorems \protect\ref{thupper-bound}
and \protect\ref{thupper-bound-2}}{Proofs of Theorems 2.1 and 2.3}}
\label{secproof-1-2}

Proofs of Theorems \ref{thupper-bound}
and \ref{thupper-bound-2} go along the same lines and exploits
three basic
statements presented here.
Lemmas \ref{lembias}~and~\ref{lemvariance}
given below
establish upper
bounds on the bias and variance of the estimator $\tilde{F}_\lambda$.
Then we\vspace*{1pt} present
Lemma \ref{lemexp-inequality} that states an
exponential inequality on the stochastic error of $\tilde{F}_\lambda
$. This
result is used for derivation bounds on the
$\epsilon$-risk.
Finally we briefly explain how the stated results
are combined in order to complete the proof of Theorems \ref{thupper-bound}
and~\ref{thupper-bound-2}.\looseness=-1

We start with the standard decomposition of the error
of estimator (\ref{eqcdfEstimator1}).
\begin{eqnarray*}
|\tilde{F}_\lambda- F(t_0)| &\leq&| \bE
\tilde{F}_\lambda
-F(t_0)| + |\tilde{F}_\lambda- \bE\tilde{F}_\lambda|
= B_\lambda(t_0;F) + |V_\lambda|,
\\
\bE|\tilde{F}_\lambda- F(t_0)|^2 &=&
B^2_\lambda(t_0;F) + \bE|V_\lambda|^2,
\end{eqnarray*}
where we have denoted
\[
B_\lambda(t_0, F) :=\biggl|\frac{1}{\pi}\int_\lambda^\infty\frac
{1}{\omega}
\Im(e^{-i\omega t_0}\hat{f}(\omega)) \,\rd\omega\biggr|,\qquad
V_\lambda:=\frac{1}{n}\sum_{j=1}^n [\xi_j(\lambda) - \bE\xi
_j(\lambda)]
\]
and
\[
\xi_j(\lambda) :=\frac{1}{\pi}\int_0^\lambda\frac{1}{\omega}
\Im\biggl\{\frac{e^{i\omega(Y_j-t_0)}}{\hat{f}_\zeta(\omega
)}\biggr\}
\,\rd\omega,\qquad
j=1,\ldots, n.
\]
%
\subsubsection{Bounds on bias and variance}
\label{subsecbounds-bias-var}
First we bound the bias of
$\tilde{F}_\lambda$.
\begin{lemma}\label{lembias}
Let $\tilde{F}_{\lambda}$ be the estimator defined in (\ref
{eqcdfEstimator1});
then for every class $\cF_\alpha(L)$ with $\alpha>-\frac{1}{2}$, $L>0$
and for any $\lambda\geq1$ one has
%
%
\begin{equation}\label{eqK-0}\qquad
\sup_{F\in\cF_\alpha(L)}
B_\lambda(F; t_0) \leq
K_0
L \lambda^{-\alpha- 1/2},\qquad
K_0:=\sqrt{2/\pi}[1+ (2\alpha+1)^{-1/2}].
\end{equation}
\end{lemma}

Now we establish an upper bound on the variance of $\tilde{F}_\lambda$.
Recall that $\omega_1$ and $c_*$
are given in (\ref{eqc-*}) and depend on
the constants $\omega_0$, $b_\zeta$ and $\tau$ appearing
in assumption~\ref{assE2}.
Define
%
%
\begin{equation}\label{eqvarphi}
w(\lambda):=\cases{
\lambda^{2\beta-1}, &\quad$\beta>1/2$,\cr
1 \vee\ln(\lambda/\omega_1), &\quad$\beta=1/2$,\cr
1, &\quad$\beta\in(0,1/2)$.}
\end{equation}
\begin{lemma}\label{lemvariance}
Let assumptions \ref{assE1}, \ref{assE2} hold and $\tilde
{F}_{\lambda}$ be
the estimator defined in (\ref{eqcdfEstimator1}).
Then there exist
constants
$K_1=K_1(\alpha,\beta,\omega_1)$ and $K_2=K_2(\beta,\omega_1)$
such that for
every $\lambda\geq1\vee\omega_1$
the following statements hold:
\begin{longlist}[(iii)]
\item[(i)] If $\alpha+\beta> 1/2$, then
\[
\operatorname{var}
\{\tilde{F}_{\lambda}\} \leq
K_1LC_\zeta c_\zeta^{-2} w(\lambda)n^{-1} + c_*n^{-1}.
\]
If $\beta>1$, then the upper bound can be made independent of
$\alpha$ and $L$
\[
\operatorname{var}
\{\tilde{F}_{\lambda}\} \leq K_2 C_\zeta
c_\zeta^{-2} \lambda^{2\beta-1} n^{-1} + c_*n^{-1}.
\]
\item[(ii)] If $\alpha+\beta=1/2$, then
\[
\operatorname{var}
\{\tilde{F}_{\lambda}\} \leq
K_1LC_\zeta c_\zeta^{-2}
w(\lambda)\sqrt{\ln(\lambda/\omega_1)} n^{-1} + c_*n^{-1}.
\]
\item[(iii)] If $\alpha+\beta<1/2$, then
%
%
\begin{equation}\label{eqzone-2}
\operatorname{var}
\{\tilde{F}_{\lambda}\} \leq
K_1 c_\zeta^{-2}
\min[
LC_\zeta\lambda^{{1/2}-\beta-\alpha},
\ln^2(\lambda/\omega_1) +\lambda^{2\beta}
]n^{-1} +
c_*n^{-1}.
\end{equation}
\end{longlist}
Explicit expressions for $K_1$ and $K_2$
are given in the proof; see \textup{(A.12)} in \cite{DGJ-suppl}.
\end{lemma}

It is worth noting that if $\beta>1$, then the upper bound on the
variance of
$\tilde{F}_\lambda$ stated in part (i) does not depend on paramaters
$\alpha$
and $L$. This is particularly important when the problem of adaptive estimation
of $F(t_0)$
is considered.

\subsubsection{An exponential inequality}
\label{subsubsecexp}
First we recall some notation.
\[
\sigma_\lambda^2 :=
c_* +
\frac{2}{\pi^2}\bE\biggl(\int_{\omega_1}^\lambda\frac{1}{\omega}
\Im\biggl\{\frac{e^{i\omega(Y_j-t_0)}}{\hat{f}_\zeta(\omega)}
\biggr\} \,\rd\omega
\biggr)^2,
\]
where
$\omega_1=\min\{\omega_0, (2b_\zeta)^{-1/\tau}\}$,
$c_*= 2\pi^{-2}[2+ (1/\tau)]^2$
and constants $\omega_0$, $b_\zeta$ and $\tau$
appear in assumption \ref{assE2}.
Define
\[
m_\lambda:= \sqrt{2c_*} + 2^{1+(\beta/2-1)_+} (\pi c_\zeta)^{-1}
[\ln(\lambda/\omega_1) + \beta^{-1} \lambda^\beta].
\]
It is easily seen that $m_\lambda\leq\bar{m}\lambda^\beta$,
$\forall
\lambda\geq1$, where
$\bar{m}$ is defined in (\ref{eqm-bar}).
%
%
We also put
\[
\bar{\sigma}^2:= c_* + C_\zeta c_\zeta^{-2}
\{ K_1L \ind(\beta\leq1) + [(K_1 L) \vee K_2]\ind(\beta
>1)\},
\]
where
constants $K_1$ and $K_2$ are given in (A.12) in \cite{DGJ-suppl}.
\begin{lemma}\label{lemexp-inequality}
Suppose that assumptions \ref{assE1} and \ref{assE2} hold; then
for any $\lambda>0$ and $z >0$ one has
%
%
\begin{equation}\label{eqV-lambda}
\bP\{|V_\lambda|\geq z\}\leq
2\exp\biggl\{-\frac{nz^2}{2\sigma_{\lambda}^2+
({2}/{3})m_\lambda z}\biggr\}.
\end{equation}
In particular, if $\alpha+\beta>1/2$, then for any $\lambda\geq
1\vee\omega_1$
and $z>0$
one has
%
%
\begin{equation}\label{eqBernstein}
\bP\{|V_\lambda| \geq z\}
\leq2 \exp\biggl\{-\frac{nz^2}{2\bar{\sigma}^2 w(\lambda)+({2}/{3})
\bar{m}\lambda^\beta z}\biggr\},
\end{equation}
where $w(\lambda)$ is given in (\ref{eqvarphi}).
\end{lemma}

\subsubsection{\texorpdfstring{Outline of the proofs of Theorems \protect\ref{thupper-bound}
and \protect\ref{thupper-bound-2}}{Outline of the proofs of Theorems 2.1 and 2.3}}
\label{subsubsecoutline}

The upper bounds on the quadratic risk stated in
Theorems \ref{thupper-bound}
and \ref{thupper-bound-2}
are immediate consequence of
Lemmas \ref{lembias} and \ref{lemvariance}.
Balancing the upper bounds on the bias and variance with respect to the
smoothing parameter $\lambda$, we come to the announced results.
Lemma~\ref{lemexp-inequality} along with Lemma \ref{lembias} are used
in order to derive upper bounds on the $\epsilon$-risk. Full technical
details are provided in the supplementary paper \cite{DGJ-suppl}.

\subsection{\texorpdfstring{Proof of Theorem \protect\ref{thadaptive}}{Proof of Theorem 2.4}}

\label{subsecadaptive00} The next preparatory lemma\vspace*{1pt}
establishes an exponential probability inequality on deviation of
$\tilde{\Sigma}_\lambda^2$ from $\Sigma_\lambda^2$.
\begin{lemma}\label{lemsigma-lambda}
Suppose that assumptions \ref{assE1} and \ref{assE2} hold.
\begin{longlist}[(ii)]
\item[(i)]
For every $\lambda\in\Lambda$
\[
\bP\{|\tilde{\Sigma}_\lambda^2 -\Sigma_\lambda^2|\geq
v_\lambda^2/2 \} \leq\frac{\epsilon}{2N}.
\]
\item[(ii)]
Let $q(\epsilon):=\sqrt{3\ln(4N\epsilon^{-1})}$; then for every
$\lambda\in\Lambda$
\[
\bP\{|V_\lambda|\geq q(\epsilon) v_\lambda n^{-1/2}\} \leq
\frac{\epsilon}{2N} .
\]
\end{longlist}
\end{lemma}

Proof of Lemma \ref{lemsigma-lambda} is given in \cite{DGJ-suppl}.

\subsubsection{\texorpdfstring{Proof of Theorem
\protect\ref{thadaptive}}{Proof of Theorem 2.4}}

Define the following events:
\begin{eqnarray*}
A(\lambda)&:=&\{|V_\lambda| \leq q(\epsilon)v_\lambda
n^{-1/2}\}
\cap\{
|\tilde{\Sigma}^2_\lambda-\Sigma_\lambda^2|\leq v^2_\lambda/2
\},
\\
A(\Lambda)&:=& \bigcap_{\lambda\in\Lambda} A(\lambda).
\end{eqnarray*}
It follows from Lemma \ref{lemsigma-lambda} and $\#(\Lambda)=N$ that
$\bP\{A(\Lambda)\}\geq1-\epsilon$.
By the triangle inequality,
%
%
\begin{equation}\label{eqsum-1-2}
|\tilde{F}_{\tilde{\lambda}}-F(t_0)| \leq|\tilde{F}_{\lambda_o} -
F(t_0)| +
|\tilde{F}_{\tilde{\lambda}} - \tilde{F}_{\lambda_o}|.
\end{equation}

By definition of $\lambda_o$ and by the fact that
$v_\lambda$ is monotone increasing with $\lambda$
we have that $v_\lambda n^{-1/2}\geq\bar{B}_\lambda$
for all $\lambda\geq\lambda_o$,
where
we have denoted\vspace*{1pt}
$\bar{B}_\lambda:=2(2/\pi)^{1/2}L\lambda^{-\alpha-1/2}$.
Therefore, on the event $A(\Lambda)$
%
%
\begin{equation}\label{eqfirst-term}
|\tilde{F}_{\lambda_o}-F(t_0)| \leq\bar{B}_{\lambda_o} +
|V_{\lambda_o}| \leq
[1+
q(\epsilon)] v_{\lambda_o} n^{-1/2}.
\end{equation}

Furthermore,
if $A(\Lambda)$ holds,
then for any pair $\lambda, \mu\in\Lambda$ satisfying
$\lambda\geq\lambda_o$ and $\mu\geq\lambda_o$ one has
$Q_\lambda\cap Q_\mu\ne\varnothing$. Indeed, by definition of
$\lambda_o$
for any $\lambda\geq\lambda_o$ one has $\bar{B}_\lambda\leq
v_\lambda/\sqrt{n}$; therefore
\[
|\tilde{F}_\lambda- F(t_0)| \leq\bar{B}_\lambda
+ q(\epsilon)v_\lambda n^{-1/2} \leq[1+q(\epsilon)]v_\lambda n^{-1/2}.
\]
In addition,
on the set $A(\Lambda)$ we have
\[
|\tilde{v}_\lambda- v_\lambda| \leq|\tilde{v}_\lambda^2
-v_\lambda^2|^{1/2}
=|\tilde{\Sigma}^2_\lambda- \Sigma^2_\lambda|^{1/2} \leq v_\lambda
/\sqrt{2}.
\]
This yields
\begin{eqnarray*}
|\tilde{F}_\lambda- F(t_0)| &\leq&\frac{\sqrt{2}}{\sqrt{2}-1}
[1+q(\epsilon)]
\tilde{v}_\lambda n^{-1/2}\\
&=&\vartheta\tilde{v}_\lambda n^{-1/2}\qquad
\forall\lambda\geq\lambda_o.
\end{eqnarray*}
Thus one has
$F(t_0)\in Q_\lambda$ and $F(t_0)\in Q_\mu$ for all $\lambda\geq
\lambda_o$
and $\mu\geq\lambda_o$; hence $Q_\mu\cap Q_\lambda\ne\varnothing$.
Then by the procedure definition, $\tilde{\lambda}\leq\lambda_o$
and $Q_{\tilde{\lambda}}\cap Q_{\lambda_o}\ne\varnothing$ on the event
$A(\Lambda)$. Therefore
%
%
\begin{eqnarray}\label{eqsecond-term}
|\tilde{F}_{\tilde{\lambda}} - \tilde{F}_{\lambda_o}|
&\leq&
\vartheta n^{-1/2} [\tilde{v}_{\tilde{\lambda}} +\tilde{v}_{\lambda_o}]
\nonumber\\
&\leq&
2\vartheta n^{-1/2}\tilde{v}_{\lambda_o}\\
&\leq&
\sqrt{2}\bigl(1+\sqrt{2}\bigr) \vartheta n^{-1/2} v_{\lambda_o}.
\nonumber
\end{eqnarray}
Here the second line follows from $\tilde{v}_{\tilde{\lambda}}\leq
\tilde{v}_{\lambda_o}$,
and the fact that
$\tilde{v}_{\lambda_o} \leq(1+2^{-1/2})v_{\lambda_o}$
on the event $A(\Lambda)$.
Combining (\ref{eqsecond-term}), (\ref{eqfirst-term}) and (\ref
{eqsum-1-2})
we obtain
that on the set $A(\Lambda)$
\[
|\tilde{F}_{\tilde{\lambda}}- F(t_0)| \leq
\biggl(\frac{3\sqrt{2}-1}{\sqrt{2}-1}\biggr)
[1+q(\epsilon)] v_{\lambda_o} n^{-1/2}.
\]
This completes the proof.

\section*{Acknowledgments}

The authors are grateful to the Associate Editor and two anonymous
referees for careful reading and useful remarks that led to
improvements in the presentation.

\begin{supplement}[id=suppA]
\stitle{Supplement to ``On deconvolution of
distribution functions''\\}
\slink[doi,text={10.1214/11-AOS907SUPP}]{10.1214/11-AOS907SUPP}
\sdatatype{.pdf}
\sfilename{aos907\_supp.pdf}
\sdescription{In the
supplementary paper \cite{DGJ-suppl} we prove results stated here and
provide additional details for the proofs appearing in
Section \ref{secproofs}. In particular, \cite{DGJ-suppl} is
partitioned in two Appendices, A and B. Appendix A contains proofs for
Section \ref{secminimax}: full technical details for
Theorems \ref{thupper-bound}, \ref{thupper-bound-2}
and \ref{thadaptive} are presented, and the proof of
Theorem \ref{thlower-1} is given. In Appendix B we prove
Theorem~\ref{thethemain1} from Section \ref{secminimax-linear}.}
\end{supplement}

%

\printaddresses

\end{document}